\documentclass[a4paper,12pt]{article}
\usepackage{amsmath,amsfonts,amssymb,amscd, array}

\def\N{{\mathbb N}}
\def\Z{{\mathbb Z}}
\def\Q{{\mathbb Q}}

\def\a{\alpha}
\def\b{\beta}
\def\l{\lambda}

\def\A{{\mathbb A}}
\def\B{{\mathbb B}}
\def\moins{\raise 1pt\hbox{{$\scriptstyle -$}}}
\def\plus{\raise 1pt\hbox{{$\scriptstyle +$}} }

\def\Sym{{ \mathfrak{ S\hspace{-0.05 em}y\hspace{-0.05 em} m}}}
\def\mfS{\mathfrak{S}}

\def\ch{\mathfrak{c\hspace{-0.07 em}h}}
\def\sh{\mathfrak{s\hspace{-0.07 em}h}}

\def\Part{{ \mathfrak{ P\hspace{-0.07 em}a\hspace{-0.05 em}r%
 \hspace{-0.05 em}t}}}

\def\nn{{\mathbf n}}

\newtheorem{example}{Example}[section]

\newtheorem{theorem}[example]{Theorem}
\newtheorem{corollary}[example]{Corollary}

\newtheorem{proposition}[example]{Proposition}

\def\Proof{\noindent {\emph Proof.}\quad }
\def\qed{\hfill $\square$}
\def\QED{\hfill $\square$}

\def\CARRE#1{\hbox{\vrule width \thickness
   \vbox to \carresize{\hrule height \thickness\vss
      \hbox to \carresize{\hss#1\hss}
   \vss\hrule height\thickness}
\unskip\vrule width \thickness}
\kern-\thickness}
\def\vsquare#1{\vbox{\CARRE{$#1$}}\kern-\thickness}

\def\smallyoung#1{%
  \ifx\carresize\undefined \newdimen\carresize \fi
  \ifx\thickness\undefined \newdimen\thickness \fi
  \carresize=14pt%
  \thickness=0.5pt%
  \vcenter{%
    \vbox{\smallskip\offinterlineskip%
      \halign{&\vsquare{##}\cr #1}}}}

\def\young#1{%
  \ifx\carresize\undefined \newdimen\carresize \fi
  \ifx\thickness\undefined \newdimen\thickness \fi
  \carresize=16pt%
  \thickness=0.5pt%
  \vcenter{%
    \vbox{\smallskip\offinterlineskip%
      \halign{&\vsquare{##}\cr #1}}}}

\def\bigyoung#1{%
  \ifx\carresize\undefined \newdimen\carresize \fi
  \ifx\thickness\undefined \newdimen\thickness \fi
  \carresize=20pt%
  \thickness=0.7pt%
  \vcenter{%
    \vbox{\smallskip\offinterlineskip%
      \halign{&\vsquare{##}\cr #1}}}}
\newdimen\unit
\def\o{$\scriptscriptstyle{{\rm o}}$}
\def\grape(#1,#2)#3{\raise#2\unit\rlap{\kern#1\unit #3}\ignorespaces}

\def\gg{{\unit=1mm
\hbox {\grape(2,1.2){'}
       \grape(1,2)\o
       \grape(2,2)\o
       \grape(3,2)\o
       \grape(1.5,1)\o
       \grape(2.5,1)\o
       \grape(2,0)\o\
}\kern 3.5 \unit}}

\def\gfill{\leaders\hbox to 1.2em{\hss\gg\hss}\hfill}
\def\frise{\medskip \centerline{\hbox to 8cm{\gfill} }\bigskip}

\date{}

\begin{document}

\title{\large\sc 
Adding $\pm 1$ to the argument of an Hall-Littlewood polynomial}

\maketitle

\frise
\begin{abstract}
Shifting by $\pm 1$ powers sums: $p_i \to p_i \pm 1$ induces a transformation
on symmetric functions that we detail in the case of Hall-Littlewood 
polynomials. By iteration, this gives a description of these
polynomials in terms of plane partitions, as well as some
generating functions. We recover in particular an identity   
of Warnaar related to Rogers-Ramanujan identities.
\end{abstract}

\frise

\hfill{\emph{D\'edi\'e \`a Xavier Viennot}

\section{Introduction}

To free analysis from infinitesimal quantities,
and quieten down the metaphysical anguish of Bishop 
Berkeley\footnote{Il s'est \'elev\'e un Docteur ennemi 
de la Science qui a d\'eclar\'e
la guerre aux Math\'ematiciens; ce Docteur monte en Chaire pour apprendre
aux fid\`eles que la G\'eom\'etrie est contraire \`a la religion; il leur dit
d'\^etre en garde contre les G\'eom\`etres, ce sont, selon lui, des gens
aveugles et indociles qui ne savent ni raisonner, ni croire;
des visionnaires qui se refusent aux choses simples et qui donnent
t\^ete baiss\'ee dans les merveilles.... , 
pr\'eface de Buffon \cite{Newton}.
Other (anonymous) defenses of infinitesimals have been published in England.
See, for example,
\emph{A defence of free-thinking in Mathematicks}, London (1735), and 
\emph{Geometry no freind to infidelity}, London (1734).},
Lagrange \cite{Lagrange} proposed to replace 
differential calculus by the study of the behaviour of a 
function\footnote{Lagrange, in harmony with D. Knuth,
 writes $fx$ and $f(x+y)$, forgetting parentheses when the argument
is a single letter. We shall follow Lagrange, when no ambiguity is to be feared.}
$fx$ under the addition of an ``increment'' $y$ to $x$:
$$   fx  \to f(x+y) \ . $$
This works perfectly well, at least under the mild proviso 
of analycity in the neighbourhood of $x$. 

One can adopt the same strategy in the realm of symmetric polynomials,
adding a letter $x$ to a set of indeterminates (we say \emph{alphabet} $A$).
However, there exists a canonical involution on symmetric functions,
and this entails that there exists in fact two versions 
of the ``addition'' of $x$.
In terms of $\l$-rings, one has to study both transformations
\begin{equation*}
  fA \to f(A+x)  \qquad ,\qquad fA \to f(A -x)  \, ,
\end{equation*}
which may look rather different when considering
explicit polynomials.

Of course, when restricting to homogeneous polynomials, one does not
lose any information by specializing $x$ to $1$, i.e. by studying instead  
\begin{equation*}
  fA \to f(A+1)  \qquad ,\qquad fA \to f(A -1)  \, .
\end{equation*}

We consider these two operations in the case of Hall-Littlewood polynomials
in Theorem \ref{th:Add1} and Theorem \ref{th:1-X}.

An interesting outcome is a short proof of an identity of Warnaar 
(\ref{th:Warnaar})
concerning a generating function of Hall-Littlewood polynomials
related to the Rogers-Ramanujan identities.

\medskip
Recall that 
the ring of symmetric polynomials $\Sym(X)$ in an infinite set of
indeterminates $X=\{x_1,x_2,\ldots\}$,
with coefficients in $\Q[[t]]$, $t$ another indeterminate,
 admits a linear basis,
the \emph{Schur functions} $S_\l X$, indexed by 
all partitions\footnote{  
We follow Macdonald's conventions \cite{Mac}.  Partitions are  weakly
decreasing sequences of non negative integers. One identifies 
two such sequences differing by adjunction of $0$'s on the right.}
$\l$. 

One can take, for a symmetric function, more general arguments
than a set of indeterminates, 
and we shall need to use, given two sets of indeterminates $X$,$Y$,
and a symmetric function $f$, 
the functions  
$f(X\pm Y)$, $f(XY)$, $f(X(1-t))$, $f(X\pm 1)$. 
Since $\Sym(X)$ is a ring of polynomials in the power sums $p_i$,
the above  symmetric functions  are induced 
from the case where $f$ is a power sum, setting
\begin{multline*}
 p_i(X\pm Y) =p_i X \pm p_i Y \, ,\  p_i(XY) = p_iX p_iY\, ,   \\  
 p_i(X(1-t))= (1-t^i)p_i X\, ,\  
 p_i(X\pm 1) = p_i X \pm 1 \, . 
\end{multline*} 

For more informations about the flexibility of 
arguments of symmetric functions,
and the use of $\l-$rings, see \cite{Cbms}.

We shall also need the generating function of  complete functions~:
$$\sigma_1 X := \prod_{x\in X} (1-x)^{-1} 
= \sum_i S_i X\, . $$

Schur functions occur naturally when decomposing the Cauchy kernel
$\sigma_1 (XY)$. The Hall-Littlewood polynomials, our present concern,
are associated
to the kernel 
$\sigma_1 (XY(1-t)):= \prod_{x\in X,y\in Y} (1-txy)(1-xy)^{-1} $.

\section{Hall-Littlewood Polynomials}

From now on, we fix a positive integer $n$. $\Part$ will be the set of
partitions of length not more than $n$, considered as elements of 
$\N^n$. 
Schur functions may be defined as determinants of complete functions,
and this allows to extend their indexation to any $v\in \Z^n$, $n$ arbitrary.
This amounts introducing the relations\footnote{
For example, $ S_{-2,2,0,0} = -S_{1,-1,0,0}= S_{1,-1,0,0,0} =0   $ 
and $S_{-3,2,1,1}= -S_{1,-2,1,1}= S_{1,0,-1,1}= -S_{1,0,0,0}= -S_1$.
}
\begin{equation} \label{Redresse1}
  S_v = - S_{\ldots, v_{i+1}-1,v_i+1,\ldots} \ ,
   \ , \ S_{v_1,\ldots,v_n}=0 \ \text{if}\ v_n<0 \ .
\end{equation}

A more powerful point of view than using straightening relations 
is to introduce 
\emph{symmetrizing operators}, which can be defined as products
of isobaric divided differences $\pi_i$, $i=1,2,\ldots$
 (operators act on their left)~:
$$ f \to f\, \pi_i :=
\frac{x_if -x_{i+1}f^{s_i}}{x_i-x_{i+1}} $$
where $s_i$  acts on functions of $X$ 
by transposition of $x_i, x_{i+1}$.

The operator $x^v\to S_v(x_1,\ldots, x_n)$, $v\in \N^n$,
 can be expressed as a product, denoted $\pi_\omega$, of
$\pi_i$. It can also be written as a summation on the symmetric group
$\mfS_n$:
$$ f\, \pi_\omega = \sum_{\sigma\in\mfS_n}  
\left(f\, \prod_{1\leq i<j\leq n} (1-x_j/x_i)^{-1}   \right)^\sigma   \, .$$
We refer to chapter 7 of \cite{Cbms} 
for some of its properties.

In particular, for any $i$, one has $\pi_i\, \pi_\omega =\pi_\omega$,
and the reordering (\ref{Redresse1})
 of the indexation of Schur functions comes from
the relation $x_{i+1}\, \pi_i=0$. 
Indeed, if $v\in \Z^n$ is such that $v_i=v_{i+1}$, and $\a,\b\in \Z$, then 
$$ x^v \left(x_i^\a x_{i+1}^\b + x_i^\b x_{i+1}^\a \right)\, x_{i+1}\, \pi_i
 = x_{i+1}\, \pi_i\, x^v \left(x_i^\a x_{i+1}^\b + x_i^\b x_{i+1}^\a \right)
=0 \, ,$$
because symmetric functions in $x_i,x_{i+1}$ commute with $\pi_i$.

Some care is needed when taking exponents or indices in $\Z^n$,
instead of $\N^n$ (this corresponds to the difference 
between using characters of the symmetric
group, or of the linear group).

We first extend the natural order on partitions to elements of 
$\mathbb{Z}^n$ by
$$ v\ge u\, \quad \text{iff}\quad \forall k> 0\,,\quad
\sum_{i=k}^{n}(v_{i}-u_{i}) \ge 0\, .$$

The operator $\pi_\omega$ commutes with multiplication by any power
of $x_1\cdots x_n$. In the case of a positive power, one has 
$(x_1\cdots x_n)^k S_\lambda = S_{\lambda +[k,\ldots,k]}$. 
But in the case of a negative power, we have also to obey
 the rule that $S_v=0$ when $v_n=0$.
The solution is to combine the symmetrization  with a truncation operator:
\begin{equation*}
  x^v \to x^v   \quad \text{if}\quad v\geq 0 \quad ,\quad 
 x^v \to 0 \quad \text{otherwise.}
\end{equation*}

Denote by $\Cup$ the operator "truncation followed by $\pi_\omega$".
Now, one has $x^v \Cup =S_v$, for all $v\in \Z^n$,
and moreover, one can compute the image of a Laurent series with
a finite number of terms of exponents $\geq 0$.

Introducing an extra indeterminate $t$, and given any  $u\in\Z^n$,
one defines the \emph{modified Hall-Littlewood polynomial} 
$Q'_u$ by
\begin{equation}
Q'_u= x^u \, \prod_{1\leq i<j\leq n}(1-tx_i/x_j)^{-1}\, \Cup \, .
\end{equation}

This is, in fact, a finite sum of Schur functions, because 
we first eliminate in the expansion of $x^u \prod(1-tx_i/x_j)^{-1}$
all the monomials which are not $\geq 0$.

The set $\{ Q'_\l:\, \l\in\Part\}$ is a basis of $\Sym$, which specializes
to the basis of Schur functions for $t=0$.
Any $Q'_u$ can be expressed in terms of the $Q'_\l$ thanks to the
following relations~:
\begin{equation}  \label{Redresse2}
 Q'_{\ldots,\a,\b+1,\ldots} + Q'_{\ldots,\b,\a+1,\ldots}
- t\, Q'_{\ldots,\b+1,\a,\ldots} -t\, Q'_{\ldots,\a+1,\b,\ldots}  =0 \, .
\end{equation}
These relations still result from $x_{i+1}\, \pi_i=0$,
because, when $v$ is such that $v_i=v_{i+1}$, then  
$$ x^v \left(x_i^\a x_{i+1}^\b + x_i^\b x_{i+1}^\a \right)
 \frac{1-tx_i/x_{i+1}}{\prod_{j<h} 1-tx_j/x_h} \, x_{i+1}\, \pi_i =0  \, ,$$
the factor on the left of $x_{i+1}$ being symmetrical in $x_i,x_{i+1}$,
and therefore commuting with $\pi_i$. 

Relations (\ref{Redresse2}), together with $Q'_{\ldots, u_n}=0$ if $u_n<0$,
suffice to express any $Q'_u$ in terms of the $ Q'_\l:\, \l\in\Part$.

The other types of Hall-Littlewood polynomials are 
\begin{eqnarray}
Q_\l X &:=& Q'_\l \bigl( X\, (1-t)  \bigr) \, ,  \\
P_\l X  &:=&  b_\l^{-1}\, Q_\l X  \, ,    
\end{eqnarray}
with $b_\l= \prod_i \prod_{j=1}^{m_i} (1-t^j) $, for 
$\l= 1^{m_1} 2^{m_2} 3^{m_3}\cdots$. 

These functions can also be defined by symmetrization, but problems
arise when taking general weights instead of only 
\emph{dominant weights} (i.e. partitions), 
see the terminal note.

We give in Corollary \ref{Plane Partitions}
a description of the functions $Q'_\l$ in term of \emph{plane partitions}.
Recall that the $Q'_\l$  
admit another combinatorial description,
this time in terms of tableaux~:
\begin{equation}
Q'_\mu = \sum_T t^{\ch T}\, S_{\sh T}   \, ,
\end{equation}
sum over all tableaux of weight $\mu$, $\sh T$ being the shape of $T$,
and $\ch T$ being the \emph{charge} of $T$ 
(cf. \cite[III.6]{Mac})\footnote{
The charge is  a \emph{rank function} on the set of 
all tableaux.  More generally, it can be defined as a function 
on words in letters $1,2,3,\ldots$,
satisfying the following relations in the case of a dominant evaluation, i.e.
when $\mu_1\geq \mu_2 \geq \mu_3 \geq \cdots$ ~:
\begin{equation}
 \ch\left(\cdots 3^{\mu_3} 2^{\mu_2} 1^{\mu_1}  \right) = 0 
  \ , \
 \ch(w\, i) =\ch(iw)+1 \ \text{if}\ i>1 \, , 
\end{equation}
and the invariance with respect to the plactic relations \cite{Lothaire}~:
\begin{equation}
  w \equiv w'  \ \Rightarrow \ \ch w = \ch w'   \ .
\end{equation}
}.

\section{Adding $1$}

Adding $\pm 1$ to the argument of a symmetric function is,
as we already said, induced by the transformation 
$$ p_i \to p_i\pm 1\, ,\, i=1,2,\ldots $$
of the powers sums $p_i$.

In terms of Schur functions, for a partition $\l\in \N^n$, this amounts 
to 
\begin{eqnarray}   \label{Pieri1}
 S_\l(X-1)  &=& \sum_{v\in \{0,1\}^n} (-1)^{|v|} S_{\l-v} X  \, , \\
 \label{Pieri2}
  S_\l(X+1)  &=& \sum_{u\in\N^n} S_{\l+u} X \, .
\end{eqnarray}

Both summations can be reduced to a summation on partitions,
erasing vertical or horizontal strips from the diagram of $\l$
\cite[I.5]{Mac}.

One can rewrite (\ref{Pieri1},\ref{Pieri2}) as 
\begin{eqnarray}   \label{Pieri3}
 S_\l(X-1)  &=& x^{\l}\, \prod_{1\leq i\leq n} (1-1/x_i)\, \pi_\omega \, , \\
 \label{Pieri4}
  S_\l(X+1)  &=& 
    x^{\l}\, \prod_{1\leq i\leq n} \frac{1}{1-1/x_i} \, \pi_\omega\, ,
\end{eqnarray}
and therefore,
\begin{eqnarray}   \label{Qmoins1}
 Q'_\l(X-1) &=& x^{\l} \frac{\prod_{1\leq i\leq n} (1-1/x_i)}{
\prod_{1\leq i<j\leq n}(1-tx_i/x_j)}  \pi_\omega 
 = \sum_{v\in \{0,1\}^n} (-1)^{|v|} Q'_{\l-v} X 
 \\
 \label{Qplus1}
 Q'_\l(X+1) &=& x^{\l}\, \frac{1}{\prod_{1\leq i\leq n} (1-1/x_i)}
\frac{1}{ \prod_{1\leq i<j\leq n}(1-tx_i/x_j)} \, \pi_\omega \,
\end{eqnarray}

The first summation is easy to reduce, using the reordering
\begin{equation}  \label{Qmoins2}
 \sum_{v} Q'_{k^m- v} = {m \brack \a} Q'_{k^{m-\a}, (k-1)^\a}   \, ,
\end{equation}
where $v$ runs over all permutations of $[1^\a,\, 0^{m-\a}]$,
and where ${m \brack \a}$ denotes the $t$-binomial  
 $(1-t^m)\cdots (1-t^{m-\a+1})/(1-t)\cdots(1-t^\a)$.

Iterating (\ref{Qmoins2}), one gets 
\begin{equation} \label{Qmoins3}
 Q'_\l(X-1)=  \sum_\mu \prod_i (-1)^{\a_i} {m_i \brack \a_i}\, Q'_\mu\, ,
\end{equation}
sum over all partitions 
$\mu= 0^{\a_1} 1^{m_1-\a_1} 1^{\a_2} 2^{m_2-\a_2} 2^{\a_3} 3^{m_3-\a_3}\cdots$
differing from $\l= 1^{m_1} 2^{m_2} 3^{m_3}\cdots$ by a vertical strip. 

The second summation is a little more complicated to transform, and 
contrary to the case of Schur functions, will not restrict to erasing
 strips.

Let us first introduce \emph{skew Hall-Littlewood polynomials}
\cite[III.5]{Mac}
$Q'_{\l/\mu}$,  by
\begin{equation}
 Q'_\l(X+Y) = \sum_{\mu\in\Part} Q'_{\l/\mu} X\, Q'_\mu Y \ , 
\end{equation}
i.e. $Q'_{\l/\mu}X$ is defined as 
the coefficient of $Q'_\mu Y$ in the expansion
of the function $Q'_{\l}$ evaluated in $X+Y$.
Note that this definition makes sense for $\l\in\Z^n$, and not only 
for $\l\in\Part$.  

For any pair of partitions $\l,\mu$, let us define
$$\nn(\l/\mu) =  \sum_i (\l^\sim_i-\mu^\sim_i)(\l^\sim_i
                                               -\mu^\sim_i-1)/2 \, ,$$
where $\l^\sim$ and $\mu^\sim$ are the partitions conjugate to 
$\l$, $\mu$ respectively. In the case where $\mu=0$, then one write
$\nn(\l)$ instead of $\nn(\l/0)$. Moreover,
$\nn(\l)= 0\l_1+ 1\l_2+2\l_3+\cdots$.

The main result of this section is the following evaluation of
$Q'_{\l/\mu} 1$.

\begin{theorem}   \label{th:Add1}
Given a partition $\l$, then
$$ Q'_{\l}(X+1) = \sum_{\mu\subseteq \l} Q'_\mu X\, \aleph(\l/\mu) \, , $$
with 
\begin{equation}  \label{Qplus2} 
 Q'_{\l/\mu} 1 = \aleph(\l/\mu):= b_\mu^{-1}\, t^{\nn(\l/\mu)}\,
 (1-t^{\nu_1-0}) (1-t^{\nu_2-1})\cdots (1-t^{\nu_r-r+1})   \ ,
\end{equation}
$\nu_1,\ldots,\nu_r$ being the parts of index 
$\mu_1,\ldots, \mu_r$, $r=\ell(\mu)$, of the partition conjugate to $\l$.
\end{theorem}

One can visualize the value of $Q'_{\l/\mu} 1=\aleph(\l/\mu) $ 
in the Cartesian plane as follows: 
color in black the rightmost box of each row of the diagram of $\mu$,
write $0,1,2,\ldots$ in the successive boxes of the same column of $\l/\mu$.
Each column with  $\a$ black boxes, and $\b$ boxes above,
gives a contribution  $t^{\b \choose 2} {\a+\b \brack \a}$ to 
$Q'_{\l/\mu}(1)$, which is the product over all columns of these 
contributions.
For example, for $\l=[4432221]$, $\mu=[2211]$, then
 $\l^\sim =[7632]$, $\nu=[6677]$, $\nn(\lambda/\mu)=13$,
$b_\mu=(1-t)^2 (1-t^2)^2$,
$$Q'_{4432221/2211} 1= t^{13} (1-t^6)(1-t^{6-1})(1-t^{7-2})(1- t^{7-3})
  (1-t)^{-2} (1-t^2)^{-2}  \, , $$
$$  \young{2\cr
           1  &3 \cr 
           0 & 2 \cr
           \blacksquare & 1\cr
         \blacksquare  &0 &2 \cr
          & \blacksquare & 1 & 1 \cr
          & \blacksquare & 0 & 0 \cr
} \ \Rightarrow \
Q'_{4432221/2211} 1 = t^3{5\brack 2}\, t^6{6\brack 2}\, 
  t^3 {3\brack 0} \, t^1 {2\brack 0}  \, .
 $$

\noindent \emph{Proof of the theorem.} \  
We shall evaluate $\sum_u Q'_{\l-u} X$ by decomposing the set 
$u\in\N^n$ into two subsets, according to whether 
$u_1=0$ or not.

Let $a,j,k$ be such that 
$$ \l_1=a=\cdots = \l_k > \l_{k+1} \ , \
  \mu_1=a=\cdots = \mu_j >\mu_{j+1} \, , $$ 
and write $\l= a^k \zeta$, $\mu= a^j \eta$.

In terms of these parameters, using (\ref{Qmoins2}), 
one wants to show that 
\begin{equation} 
  Q'_{\l/\mu} 1= t^{{k-j \choose 2}} {k\brack j} 
    Q'_{[a^{k-j} \zeta]/\eta}(1)\, .
\end{equation} 

The sub-sum $\sum_{u:\, u_1=0} Q'_{\l-u} X$  
is equal to 
$$   Q'_a X \odot Q'_{a^{k-1}\zeta}(X+1) \, ,  $$
defining $ Q'_a \odot Q'_\nu$ to be the concatenation $Q'_{a,\nu}$,
and extending by linearity.

By induction on $\ell(\l)$, the coefficient of 
$Q'_\mu X$ in $   Q'_a X \odot Q'_{a^{k-1}\zeta}(X+1)$ is equal to
$$ t^{{k-j \choose 2}} {k-1\brack j-1} Q'_{[a^{k-j} \zeta]/\eta} 1 \ . $$ 

The terms $Q'_{\l -u}X$, $u_1\geq 1$, can be rewritten 
$$  \sum_{u\in\N^n} Q_{\l^- -u} X 
Q'_{\l^-}(X+1) = t^{k-1} Q'_{a^{k-1},a-1,\zeta} (X+1)  \, , $$
with $\l^- =[a\moins 1,\l_2,\l_3,\ldots]$.

By induction on $|\l|$, the coefficient of $Q'_\mu X$ in this function
is equal to
$$  
 t^{k-1} t^{{k-1-j \choose 2}}  {k-1\brack j} Q'_{[a^{k-j} \zeta]/\eta} 1\ .
$$

The identity 
$$  {k-1\brack j-1} + t^{k-1 +{k-1-j \choose 2}-{k-j \choose 2}}
  {k-1\brack j} = {k-1\brack j-1} + t^j {k-1\brack j} = {k\brack j} $$
allows to sum up the two contributions and finishes the proof.   \qed

\def\bs{\blacksquare}
For example, $Q'_{221}(X+1)$ is obtained by the following enumeration~:
$$ t^4\, \smallyoung{2\cr 1 &1\cr 0&0\cr} + 
t^2 {3 \brack 1}\,  \smallyoung{1\cr 0 &1\cr \bs &0\cr}  
+t {2 \brack 1}\,  \smallyoung{1\cr 0 &0\cr  &\bs\cr}
+t{3 \brack 2}\,  \smallyoung{0\cr \bs &1\cr \bs &0\cr}
+{2 \brack 1}^2 \,  \smallyoung{0\cr \bs &0\cr &\bs\cr}
$$
$$
+t{3 \brack 3}\,  \smallyoung{\bs\cr \bs &1\cr \bs &0\cr}
+ {2\brack 2} \, \smallyoung{0\cr &\bs \cr &\bs \cr}
+ {2\brack 2}{2\brack 1}\, \smallyoung{\bs\cr \bs & 0\cr &\bs\cr}
+ {1\brack 1}{2\brack 2}\, \smallyoung{\bs\cr &\bs\cr &\bs\cr}
 $$
\begin{multline*}
 = t^4 Q'_{0} +t^2(1+t+t^2) Q'_1 +t(1+t)Q'_2 + t(1+t+t^2) Q'_{11}
                                                      +(1+t)^2 Q'_{21} \\
 +t Q'_{111} + Q'_{22} + (1+t) Q'_{211} + Q'_{221}  \, .
\end{multline*}

Recall that a \emph{plane partition}  of shape a partition $\l$
is a filling of the diagram of $\l$ with positive integers such 
that numbers weakly increase when faring towards the origin. In other words,
the successive domains occupied by the letters $1,2,,3,\ldots$
are skew partitions $\l/\l^1$, $\l^1/\l^2$, $\l^2/\l^3,\, \ldots$.
Define the \emph{weight} of a plane partition $\boxplus$ to be the 
product
$$ \aleph_x\left(\boxplus \right) =
  \aleph(\l/\l^1) x_1^{|\l/\l^1|}\,  \aleph(\l^1/\l^2) x_2^{|\l^1/\l^2|}\,
 \aleph(\l^2/\l^3) x_3^{|\l^2/\l^3|}\,  \cdots \ . $$

Iteration of Theorem \ref{th:Add1} leads to the following 
description of $Q'_\l$.
\begin{corollary}   \label{Plane Partitions}
Let $\l$ be a partition, $n$ be  a positive integer. Then
$$ Q'_\l(x_1+\cdots +x_n) = \sum_{\boxplus} \aleph_x\left(\boxplus \right)
 \, ,$$
sum over all plane partitions in $1,\ldots,n$ of shape $\l$.
\end{corollary}

For example, for $\l=[21]$, one has
\begin{eqnarray*}
 Q'_{21} &=& \sum_i t\, \young{i\cr i & i \cr}
 +\sum_{i<j} (1\plus t)\, \young{i\cr j & i \cr} +t\, \young{j\cr j & i \cr}
+ \young{i\cr j & j \cr} \\
 & & \hspace*{60 pt}  + \sum_{i<j<k} \young{i\cr k & j \cr}
 +(1\plus t)\, \young{j\cr k & i \cr}  \\
 &=& \sum_i t\, x_i^3 +\sum_{i<j} (1\plus t)\, x_i^2 x_j + t\, x_ix_j^2 
 + x_i x_j^2 +\sum_{i<j<k} (2+t) x_ix_jx_k \, . 
\end{eqnarray*}

On the other hand, the interpretation in term of tableaux and charge reads 
\begin{eqnarray*}
 Q'_{21} &=&  S_{21} +t S_3  \\
 &=&  \sum_{i<j} \young{j\cr i & i \cr} +\young{j\cr i & j \cr}
 + \sum_{i< j< k}\,  \young{j\cr i & k \cr} 
 + \young{k\cr i & j \cr} + \sum_{i\leq j \leq k} 
t\ \young{ i & j &k \cr} \, .
\end{eqnarray*}

\section{Argument $1 - X$}  

Instead of describing $Q'_\l(X-1)$, we prefer taking 
$Q'_\l(1-X)$. Addition and subtraction of  alphabets 
involve rather different properties of symmetric functions.
In the case of HL-polynomials, with some hypothesis on $\lambda$,
we shall find a connection between $Q'_\l(1-X)$ and resultants.  

Recall that, given two finite alphabets $\A, \B$, of respective
cardinalities $\a,\b$, then the \emph{resultant} 
$\prod_{a\in\A,b\in \B}(a-b)$
is equal to the Schur function $S_{\b^\a}(\A-\B)$.
More generally, the resultant appears as a factor of
some Schur functions,
thanks to a proposition due to Berele and Regev\cite[1.4.3]{Cbms}~:

\begin{proposition}   
 Let $\A$, $\B$, be of cardinalities $\a,\b$, $p\in \N$,
 $\zeta\in \N^p$, $\nu\in \N^\a$. Then,
writing $[\b+\nu_1,\ldots, \b+\nu_\a, \zeta_1,\ldots, \zeta_p]=
  [\b^\a+\nu,\zeta]$, one has 
\begin{equation}  \label{BereleRegev}
 S_{\b^\a+\nu,\, \zeta}(\A-\B)=
 S_\zeta(-\B)\,     S_\nu(\A)\, \prod_{a\in\A,b\in \B}(a-b)\ .
\end{equation} 
Moreover,
\begin{equation} \label{BereleRegev2}
\nu\in\Part,\,  \nu\supseteq (\b+1)^{\a+1}\
\Rightarrow S_\nu(\A-\B)=0 \, .
\end{equation}
\end{proposition}


Pictorially, the relation is

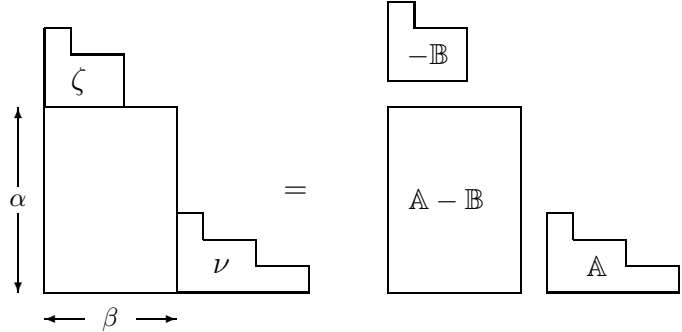
\setlength{\unitlength}{1pt}

\centerline{
\begin{picture}(190,140)(0,-20)
\put(0,0){\framebox(50,70){}}
\put(0,70){\line(0,1){30}}
\put(0,100){\line(1,0){10}}
\put(10,100){\line(0,-1){10}}
\put(10,90){\line(1,0){20}}
\put(30,90){\line(0,-1){20}}
\put(50,0){\line(1,0){50}}
\put(100,0){\line(0,1){10}}
\put(80,10){\line(1,0){20}}
\put(80,10){\line(0,1){10}}
\put(60,20){\line(1,0){20}}
\put(60,20){\line(0,1){10}}
\put(60,30){\line(-1,0){10}}
\put(69,10){\makebox(0,0){\hbox{$\nu$ }}}
\put(15,80){\makebox(0,0){\hbox{$\zeta$ }}}
\put(15,-10){\vector(-1,0){15}}
\put(35,-10){\vector(1,0){15}}
\put(25,-10){\makebox(0,0){\hbox{\footnotesize $\beta$}}}
\put(-10,30){\vector(0,-1){30}}
\put(-10,40){\vector(0,1){30}}
\put(-10,35){\makebox(0,0){\hbox{\footnotesize $\alpha$}}}
%
\put(130,0){\framebox(50,70){}}
\put(152,35){\makebox(0,0){\hbox{\small $\A-\B$}}}
\put(130,80){\line(1,0){30}}
\put(130,80){\line(0,1){30}}
\put(130,110){\line(1,0){10}}
\put(140,110){\line(0,-1){10}}
\put(140,100){\line(1,0){20}}
\put(160,100){\line(0,-1){20}}
\put(145,90){\makebox(0,0){\hbox{\small $-\B$}}}
\put(95,38){\makebox(0,0){\hbox{=}}}
 \put(190,0){\line(0,1){30}}
 \put(190,0){\line(1,0){50}}
 \put(240,0){\line(0,1){10}}
 \put(220,10){\line(1,0){20}}
 \put(220,10){\line(0,1){10}}
 \put(200,20){\line(1,0){20}}
 \put(200,20){\line(0,1){10}}
 \put(200,30){\line(-1,0){10}}
 \put(209,10){\makebox(0,0){\hbox{\small $\A$}}}
\end{picture}
}

The next theorem shows that $Q'_{n^k}(1-X)$ is also a resultant,
and more generally, factors out a resultant from $Q'_\l(1-X)$
(it is more convenient to take a power of $t$ instead of $1$,
to simplify exponents).

\begin{theorem}  \label{th:1-X}
Let $n,r\in \N$, $X$ be an alphabet of cardinality $n$,
$\l$ be a partition that is written, with $\zeta$ such that 
$\zeta_1<n$, 
$$  \l= [n+\nu_1,\ldots, n+\nu_k, \zeta_1,\zeta_2,\ldots  ]
= [n^k+\nu, \zeta] \, . $$
Then 
\begin{equation}   \label{1moinsX1}
 Q'_\l(t^r-X) = t^{\nn(\nu)+r|\nu|}\,
 \prod_{i=r}^{k+r-1} \prod_{j=1}^n (t^i-x_j)\, Q'_\zeta(t^{k+r}-X) \, .
\end{equation}
\end{theorem}

\Proof  \ $Q'_\l(1-X) = \sum_{\mu\subseteq \l} Q'_\mu X\, Q'_{\l/\mu} t^r$.
The terms $Q'_\mu X$ vanish if $\mu_1 >n$. 
Eq. \ref{Qplus2} shows that $Q'_{\l/\mu} t^r$, when $\mu\leq n$,
 is equal to
$$ Q'_{[n^k,\zeta]/\mu}t^r\, t^{\nn(\nu) +r|\nu|} \, .$$
Thus we need only treat the case where $\nu=0$, and we shall do it
by induction on $k$.

The function $Q'_{[n^{k-1},\zeta]/\mu}( t^{r+1}-X)$ can be written as 
as sum of Schur functions , enumerating all tableaux of
commutative evaluation 
$$ 2^n\cdots k^n (k\plus 1)^{\zeta_1} (k\plus 2)^{\zeta_2} \cdots \, ,$$ 
of shape not containing $[n+1,n+1]$
(otherwise the Schur function vanishes, according to (\ref{BereleRegev2})).
Given any such tableau $T$, of shape $\rho$, 
given any horizontal strip $\rho/\xi$, 
then there exists a unique pair $(u,T)$ such that 
$T\equiv u T'$, with $T'$ of shape $\xi$ and $u$ a tableau of row shape.
If $\xi_1\leq n$, then $T''= T'\, 1^n u$ is a tableau of charge
$\ch T +\ell(u)$, which contributes to
$Q'_{n^k,\zeta}$, and all tableaux of weight 
$[n^k,\zeta]$, of shape not containing 
$[n\plus 1,n\plus 1]$ are obtained in this way.

The contribution of $T$ to $Q'_{[n^{k-1},\zeta]/\mu}( t^{r+1}-X)$
is, writing $R(y,X)$ for $\prod_{1\leq i\leq n} (y-x_j)$,  
$$ t^{\ch T} S_\rho(t^{r+1}-X) = t^{\ch T} 
S_{\rho_2,\rho_3,\ldots}(-X)\,  R(t^{r+1},X)\, t^{(r+1)(\rho_1-n)} $$
thanks to the factorization (\ref{BereleRegev}). 
Each of its successor $T''= T'1^n u$ contributes to 
$Q'_[n^k,\zeta](t^r-X)$ as
$$ t^{\ch T''} S_{n+\ell(u), \xi}(t^r-X) = 
 t^{\ch T} t^{\ell(u)} S_\xi(-X) R(t^r,X) t^{r\ell(u)} \, . $$ 
Summing over all the tableaux  $T''$ whose predecessor is $T$, one gets
the contribution\footnote{Indeed,  when $y$ a single letter,
$ S_\rho(y-X)$
is equal to the sum $\sum_\xi y^{|\rho/\xi|} S_\xi(-X)$.}
$$  t^{\ch T}\, R(t^r,X)\,  S_\rho(t^{r+1}-X) \, ,$$
and this shows that the contribution of $T$ 
has been multiplied by a factor independent of $T$.

Summing over all tableaux $T$ of weight $[n^{k-1},\zeta]$ gives the
equality 
$$ Q'_{n^k,\zeta}(t^r-X) = R(t^r,X) Q'_{n^{k-1},\zeta}(t^{r+1}-X)  $$
and finishes the proof.    \qed

For example, for $n=2$, let us illustrate, on a single tableau,
 the inductive step
from $Q'_{222 11}(t-X)$ to $Q'_{2222 11}(1-X)$.
The tableau $\smallyoung{6\cr 5\cr 3 &4\cr 2&2 &3&4\cr}$
has charge $5$ and 
gives the contribution 
$$ t^5 S_{4211}(t-X)= t^5S_{211}(\moins X)\, S_4(t-X) = 
t^5S_{211}(\moins X)\, R(t,X)\, t^2 \, .$$

Its different factorizations, the issueing new tableaux and 
their contributions are 

$$ \begin{array}{lccr}
\smallyoung{2&3\cr}\ \smallyoung{6\cr 5\cr 4 &4\cr 2&3\cr} 
& \to &
\smallyoung{6\cr 5\cr 4 &4\cr 2&3\cr 1&1&2&3\cr} 
& 
 t^7 S_{2211}(\moins X)\, R(1,X) \\
\smallyoung{2&3&6\cr}\ \smallyoung{ 5\cr 4 &4\cr 2&3\cr}
& \to &
\smallyoung{ 5\cr 4 &4\cr 2&3\cr 1&1&2&3&6\cr}
&
 t^8 S_{221}(\moins X)\, R(1,X) \\
\smallyoung{2&3&4\cr}\ \smallyoung{6\cr 5\cr 4\cr 2&3\cr}
& \to &
\smallyoung{6\cr 5\cr 4\cr 2&3\cr 1&1&2&3&4\cr}
&
 t^8 S_{2111}(\moins X)\, R(1,X)\\
\smallyoung{2&3&4&6\cr}\ \smallyoung{5\cr 4\cr 2&3\cr}
& \to &
\smallyoung{5\cr 4\cr 2&3\cr 1&1&2&3&4&6\cr}
&
 t^9 S_{211}(\moins X)\, R(1,X) \, . 
\end{array}
$$

The sum of these contributions is 
\begin{multline*} t^7R(1,X)\, 
\bigl(S_{2211}(\moins X)+t\, S_{221}(\moins X)+t\, 
S_{2111}(\moins X)+t^2\,S_{211}(\moins X)\bigr) \\
= t^7R(1,X)\, S_{2211}(t- X)= t^7R(1,X) R(t,X) S_{211}(t- X) \, .
\end{multline*}

\bigskip
In the special case $X=\{x\}$ of cardinality $1$, and $r=0$, one recovers
that \cite[p.226]{Mac}
\begin{equation}  \label{1moinsX2}
 Q'_\l(1-x) = Q_\l\left(\frac{1-x}{1-t}  \right) =
t^{\nn(\l)} (1-x)(1-xt^{-1})\cdots (1-x t^{1-\ell(\l)}) \, ,
\end{equation}
identity which allows to describe the \emph{principal specializations} 
$Q'_\l(1-t^N)= Q_\l((1-t^N)/(1-t))$ when taking $y=t^N$. 

In the case of an alphabet of cardinality $2$, $X=\{x_1,x_2\}$,
the remaining partition $\zeta$ is such that 
$\zeta=1^\b$ for some $\b\in\N$, and therefore, the factor 
$Q'_\zeta(t^j -X)$ 
is equal to 
$$Q_{1^\b}\left( \frac{t^j -X}{1-t} \right) =
  b_{1^\b}\, S_{1^\b}\left( \frac{t^j -x_1-x_2}{1-t} \right) \, .$$ 
In total, 
\begin{equation}  \label{1moinsX3}
Q'_{2^k+\nu,\, 1^\b}(1-x_1-x_2) = t^{\nn(\nu)}\,
 (1\moins t)\cdots (1\moins t^\b)\, 
 \prod_{i=0}^{k-1} (t^i\moins x_1)(t^i\moins x_2)\, S_{1^\b}
   \left( \frac{t^k \moins x_1\moins x_2}{1-t} \right) \, .
\end{equation}

\section{Generating Functions}

The results of the preceding sections may be used to compute 
generating series of Hall-Littlewood polynomials.
For example, given the series 
$\sum_{\mu\in\Part} c_\mu P_\mu(X)$,
with some arbitrary coefficients $c_\mu$, suppose that 
one wants to evaluate the product
$$ \sigma_1(-X)\, \sum c_\mu P_\mu X =
 \prod_{x\in X} (1-x)\, \sum c_\mu P_\mu X \ .$$
To do so, one first extends the family $c_\mu$ to a family 
$c_v:\, v\in \Z^n$ by imposing the relations (\ref{Redresse2}). 
It amounts introducing a second alphabet $Y$, and putting $c_\mu=Q_\mu Y$.
Since $\sigma_1(XY(1\moins t)) = \sum Q_\mu Y P_\mu X$, then
\begin{multline*}
  \sigma_1(-X)\, \sum_\mu c_\mu P_\mu(X) =
   \sigma_1\bigl(-X+XY(1\moins t)\bigr) 
 = \prod_{x\in X} (1\moins x) \prod_{x\in X,y\in Y}\frac{1-txy}{1-xy}   \\
= \sigma_1\left(X(Y-\frac{1}{1-t})(1-t)   \right) =
     \sum_{\l\in\Part} P_\l X\, Q_\l\left(Y-\frac{1}{1\moins t}\right)  \, .  
\end{multline*}

Knowing the expansion of $Q_\l(Y-1/(1\moins t))$, or equivalently,
of $Q'_\l(Y'-1)$, with $Y'=Y(1\moins t)$, given in (\ref{Qmoins1}), 
one concludes
\begin{equation}  \label{Prod-X}
  \sigma_1(-X)\, \sum_\mu c_\mu P_\mu X =
  \sum_\l \sum_{v\in \{0,1\}^n} (\moins 1)^{|v|} c_{\l-v} P_\l X   \, .
\end{equation}

Let us give another similar example.
\begin{proposition}  Given two alphabets $X,Y$, then 
\begin{equation}    \label{sigmaXY}
\sigma_1\left(X+XY(1\moins t)\right) = 
\prod_{x\in X} \frac{1}{1-x} \prod_{y\in Y} \frac{1-txy}{1-xy} 
= \sum_{\l\in\Part} \sum_{\mu\subseteq \l} 
 P_\l X P_\mu Y\, b_\mu Q'_{\l/\mu} 1   \, ,
\end{equation}
the value of $b_\mu Q'_{\l/\mu} 1$ being given in  (\ref{Qplus2}).
\end{proposition}

\Proof \ One writes $X+XY(1\moins t)= X(Y+(1\moins t)^{-1})(1\moins t)$,
and therefore
$$ \sigma_1\left(X+XY(1\moins t)\right) = 
  \sum_\l P_\l X Q_\l(Y+(1\moins t)^{-1}) =
\sum_{\l,\mu}  P_\l X P_\mu b_\mu Y Q_{\l/\mu}((1\moins t)^{-1})  \, ,$$
which is the required formula.   \qed

For example, the coefficient of $P_{42}X$, in terms of the $P_\mu=P_\mu Y$,
 is
\begin{multline*}
 t^2P_{0}+ t(1\moins t^2)P_1 + (1\moins t^2)P_2 
                             +t(1\moins t)(1\moins t^2)P_{11}\\
+ (1\moins t)P_3 + (1\moins t)(1\moins t^2)P_{21}+ (1\moins t)P_4 +
(1\moins t)^2P_{31} + (1\moins t)(1\moins t^2)P_{22}   \\
+(1\moins t)^2P_{41}+ (1\moins t)^2P_{32} +(1\moins t)^2P_{42}\, .
\end{multline*}

Using the explicit values (\ref{1moinsX3}), one obtains as a corollary
for \\ $X= (t-x_1-x_2)(1-t)^{-1}$, $Y=(1-y)(1-t)^{-1}$, the expansion of~:
\begin{multline*}
\sigma_1\left( \frac{(t-x_1-x_2)(2-y)}{1-t} \right) =     
  \sigma_1\left(2 \frac{t}{1-t}+ \frac{y(x_1+x_2)}{1-t}
 -\frac{yt}{1-t} -2\frac{x_1+x_2}{1-t}  \right)   \\
 = \prod_{i=0}^\infty \frac{1}{(1-t^{i+1})^2}    
    \frac{1}{(1-t^iyx_1)(1-t^i yx_2)}  
  (1-yt^{i+1}) \bigl((1-t^ix_1)(1-t^ix_2)   \bigr)^2  \, .
\end{multline*}


Warnaar evaluates a similar generating series,
with simpler coefficients. 
Let, for any pair of partitions,
$$ \theta(\l,\mu) = t^{\nn(\l/\mu)-|\mu|} \, . $$

Then Warnaar's identity \cite[Th. 1.1]{Warnaar} is the following.

\begin{theorem}  \label{th:Warnaar}
For any pair of alphabets $X,Y$, one has 
\begin{equation}  \label{Warnaar1}
\sigma_1\left(X+Y+(\frac{1}{t}-1)XY   \right) =
\sum_{\l,\mu\in \Part} \theta(\l,\mu) P_\l X P_\mu Y  \, .
\end{equation}

\end{theorem}
\Proof  Writing the left-hand side 
$\sigma_1 X\, \sigma_1\left(Y (1+(t^{-1}-1)X)\right) $, and using that 
$\sigma_1 XY((1\moins t)= \sum Q_\mu X P_\mu Y$, one rewrites 
the identity to prove as 
$$ 
Q_\l\left(\frac{1}{1-t} +t^{-1} X  \right) 
 \stackrel{ ?}{=} 
  \frac{1}{\sigma_1 X} \sum_\mu \theta(\l,\mu) P_\mu X \, , $$ 
or 
\begin{multline}  \label{Warnaar3}
 \sum_\mu t^{-|\mu|} Q_\mu X\, Q_{\l/\mu}\left(\frac{1}{1-t}  \right)=
   \sum_\mu t^{-|\mu|} Q_\mu X\,  Q'_{\l/\mu} 1  \\
 \stackrel{?}{=} \sigma_1(\moins X)\, \sum_\mu \theta(\l,\mu) P_\mu X \, .
\end{multline}

Thanks to (\ref{Prod-X}), the right-hand side can be written 
$$ \sum_{v\in \{0,1\}^n} \theta(\l,\mu-v) P_\mu X =
 \sum_{v\in \{0,1\}^n} \theta(\l,\mu-v)   b_\mu^{-1} Q_\mu X \, .
$$
Therefore, (\ref{Warnaar3}) is a consequence of (\ref{Qplus2}) and
is in fact, equivalent to (\ref{sigmaXY}).      \QED

Let us mention that Warnaar used 
 another expression of $\theta(\l,\mu)$, putting~:
\begin{equation}    \label{DefTheta}
 \theta(\l,\mu) = t^{\nn(\l)+\nn(\mu) -(\l^\sim,\, \mu^\sim)}  \ ,
\end{equation}
where $(\l^\sim,\mu^\sim)= \sum \l_i^\sim \mu_i^\sim$.

\section{Scalar Product}

\bigskip
The combinatorics of Hall-Littlewood polynomials fundamentaly reduces 
to the  fact that $\{ P_\lambda\}$ is an orthogonal
basis of the ring of polynomials, and that it can be extended to a family
satisfying  
 the \emph{straightening relations} (\ref{Redresse2}).
Given these two ingredients, we can forget Hall and Littlewood altogether.

Thus, let us consider
the ring of Laurent series in $x_1,\ldots, x_n$, with a finite number of terms
with exponent $\geq [0,\ldots,0]$,
modulo the ideal generated by the relations 
\begin{equation}   \label{Redresse3}
 (x^v+ x^{vs_i})(x_{i+1}-tx_i)\simeq 0\, ,\, i=1,\ldots,\, n\moins 1,
   \, v\in \Z^n,\,  
 x^v\simeq 0 \ \text{if}\ v_n<0 \, .
\end{equation} 

Any element of this ring can be written uniquely as a linear combination
of dominant monomials $x^\l:\, \l_1\geq \cdots \geq \l_n\geq 0$.

We define a scalar product by its restriction to dominant monomials~: 
\begin{equation}
  (( x^\l,x^\mu))= b_\l\, \delta_{\l,\mu}\, ,\, \l,\mu\in \Part\, .
\end{equation}

We may remark that the scalar product on Laurent polynomials
used in \cite{tKey}
\begin{equation}  \label{Scalart}
 (f,g)_t := CT\left(
     f(x_1,\ldots, x_n)\,  g\left(\frac{1}{x_n},\ldots, \frac{1}{x_1}\right)
         \prod_{1\leq i<j\leq n} \frac{1-x_i/x_j}{1-tx_i/x_j} \right)\, , 
\end{equation}
where $CT$ means ``constant term'',
is such that 
$$ (Q_\l , x^\mu) = b_\l\, \delta_{\l,\mu} \, .$$
This indicates that the present study is related to more general 
constructions about non-symmetric Hall-Littlewood polynomials.

In the present set-up, to recover Warnaar's generating function, 
we essentially need to interpret the function 
$\theta(\l,\mu)$, $\l,\mu\in \Part$, as a scalar product.
To do so, let us also write $\theta(x^\l,x^\mu)$ for $\theta(\l,\mu)$, and 
extend the definition of $\theta$ to all monomials by linearity
using  relations (\ref{Redresse3}). 

\begin{proposition}
For any pair $\l,\mu\in \Part$, then 
\begin{equation}  \label{theta2}
 \left(\!\left(\,  \frac{x^\l t^{-|\l|}}{\prod_{i=1}^n 1 -t/x_i}\, ,\, 
  \frac{x^\mu}{\prod_{i=1}^n 1 -1/x_i}\, \right)\!\right) = \theta(\l,\mu)  \, .
\end{equation}

\end{proposition}

\Proof Multiplying\footnote{
Multiplying by $x^{-u}$, $u\in\N^n$, is compatible with relation 
(\ref{Redresse3}): $x^v=0$ if $v_n<0$.
}
 $x^\mu$ by $\prod_{i=1}^n 1 -1/x_i$, 
 let us prove the equivalent statement that 
\begin{equation}   \label{theta3}
 ((\,  \frac{x^\l t^{-|\l|}}{\prod_{i=1}^n 1 -t/x_i}\, ,\,
  x^\mu )) = \sum_{v\in \{0,1\}^n} (\moins 1)^{|v|} \theta(\l,\mu-v)  \, .
\end{equation}
Let $a,k$ be such that $\mu_1=a=\cdots =\mu_k>\mu_{k+1}$, 
let $r=\l^\sim_a$. Then, for any $j:\, 0\leq j\leq k$, one has
$$ \theta(\l,\mu -[1^j,0^{n-j}])= \theta(\l,\mu) 
   t^{-0-1-\cdots-(j-1)} t^{jr} \, .$$

Summing over all vectors $u\in \{0,1\}^n$, $u_i=0$ for $i>k$,
and taking into account reordering, one gets
\begin{multline*}
 \sum_u (\moins 1)^{|u|} \theta(\l,\mu-u)= 
 \theta(\l,\mu) \sum_{j=0}^k (\moins 1)^j 
             {k\brack j} t^{-{j\choose 2}} t^{jr}  \\
= \theta(\l,\mu) (1-t^{r-0})(1-t^{r-1})\cdots (1-t^{r-k+1}) \, .
\end{multline*} 

By induction $[\mu_{k+1},\mu_{k+2},\ldots] \to 
\mu=[a^k,\mu_{k+1},\ldots]$, this proves that 
\begin{equation}  \label{theta4}
 \sum_{v\in \{0,1\}^n} (\moins 1)^{|v|} \theta(\l,\mu-v)=
\theta(\l,\mu) (1-t^{\nu_1})(1-t^{\nu_2-1})(1-t^{\nu_3-2})\cdots  \, ,
\end{equation}
$\nu_1,\nu_2,\nu_3,\ldots$ being the parts of $\l^\sim$ of index $\mu_1,
\mu_2,\mu_3,\ldots$.  
Notice that the product is null if $\mu \not\subseteq \l$. 

Using (\ref{Qplus2}), one rewrites the right hand side of (\ref{theta4})
as
\begin{equation}  \label{theta5}
  \theta(\l,\mu) t^{-\nn(\l/\mu)}\, b_\mu\, Q'_{\l/\mu}\, 1  \, .
\end{equation}

On the other hand, 
$$ Q'_\l(1 +t^{-1}X) = \sum_{u\in\N^n} Q'_{\l-u}(t^{-1}X) 
       =\sum t^{-|\mu|} Q'_\mu X\, Q'_{\l/\mu} 1 \, .$$
This implies that 
$$ t^{-|\l|} \prod_{i=1}^n (1 -t/x_i)^{-1} \simeq \sum_{u\in\N^n}
x^{\l-u} t^{|u|-|\l|}  $$
be congruent to $\sum_\mu t^{-|\mu|} x^\mu Q'_{\l/\mu} 1$ .
Therefore,
$$
 ((\,  \frac{x^\l t^{-|\l|}}{\prod_{i=1}^n 1 -t/x_i}\, ,\,
  x^\mu )) =  t^{-|\mu|} b_\mu Q'_{\l/\mu} 1 $$
and this is the required identity (\ref{theta3}), using that 
$\theta(\l,\mu) = t^{\nn(\l/\mu)-|\mu|}$.   \qed

\section{Note}
The functions $Q'_\l$ have an interpretation in terms of the cohomology
of  flag manifolds \cite{DLT},
\cite[III.7]{Mac}. They describe graded multiplicities
of representations of the symmetric group. The functions $Q_\l$ also
have an interpretation, as Euler-Poincar\'e characteristic of line bundles
over the flag manifold \cite{EP}. In combinatorial terms, this gives
the following definition.

 Let $\pi_\omega$ be the symmetrizing operator defined before.
Given a partition $\l\in\N^n$, let $m_0$ be the multiplicity of the part $0$
(=$n-\ell(\l)$).  Then \cite[III.2]{Mac}
\begin{equation} \label{DefQ}
  Q_\l(X) = \frac{ (1-t)^n}{(1-t)\cdots (1-t^{m_0})}\, 
x^\l\, \prod_{i<j} (1-t x_j/x_i)\   \pi_\omega \, .
\end{equation}

However, the fact that the normalization factor depends on $\l$ prevents 
us from using (\ref{DefQ}), which is equivalent to
formula (2.14) of \cite{Mac}, as a definition for $Q_v$ when $v$ is not
a partition (contrary to what is stated p. 214 of \cite{Mac}).

Indeed, let $n=2$. Then 
$$ x^{02} (1-tx_2/x_1)\pi_\omega = 
  \frac{ x^{02} (1-tx_2/x_1)}{ 1-x_2/x_1} + 
  \frac{ x^{20} (1-tx_1/x_2)}{ 1-x_1/x_2} =
 t(x^{20}+x^{11}+x^{02}) - x^{11}  $$

On the other hand, the relations (\ref{Redresse2}) impose\footnote{
Macdonald \cite[p. 214]{Mac} writes
$Q_{16}= tQ_{61}+ (t^2-1) Q_{52}+ (t^3-t)Q_{43}$,
$Q_{15}= t Q_{51} + (t^2-1) Q_{42} +(t^2-t) Q_{33}$.
Read $t^m-t^{m-1}$ instead of $t^{m+1}-t^m$.}
$$ Q'_{02} =  t\, Q'_{20} +(t-1)Q'_{11}  $$
and the same relation must be valid for $Q_{02}$, resulting
from the transformation $X\to X(1-t)$.
But 
$$ t Q_{20}(X) +(t-1)Q_{11}(X) =
 (t-t^2)\, (x^{20}+x^{11}+x^{02}) 
 +(t-1)(1-t+t^3) x^{11}  $$
is not proportional to the image of $x^{02}$.

The confusing fact is that 
$$ ( x^{02} -t x^{20}+(1-t)x^{11})\, (1-tx_2/x_1)\, \pi_\omega =0 \ , $$
that is, the images of monomials under the operator
$(1-tx_2/x_1)\pi_\omega$ satisfy the same relations as the
functions $Q_v$. However, because of the different normalization factors,
this does not imply that the image of $x^{02}$ be proportional to
$Q_{02}$.


\bigskip
\begin{center}
{\large Alain Lascoux}\\
 CNRS, IGM, Universit\'e Paris-Est \\
  77454 Marne-la-Vall\'ee Cedex, France\\
  Email: Alain.Lascoux@univ-mlv.fr\\
\end{center}


\begin{thebibliography}{11}

\bibitem{DLT} {\sc J. D\'esarm\'enien, B. Leclerc and J.-Y. Thibon}.
\emph{Hall-Littlewood Functions and Kostka-Foulkes Polynomials 
in Representation Theory}, S\'eminaire Lotharingien de Combinatoire
{\bf 32} (1994).    


\bibitem{tKey} {\sc  F. Descouens, A. Lascoux}.
\emph{Non-symmetric Hall-Littlewood polynomials},
S\'eminaire Lotharingien de Combinatoire, SLC {\bf 54} B54Ar (2006). 

\bibitem{Lagrange} {\sc J.-L. Lagrange}.
\emph{ 
Th\'eorie des fonctions analytiques contenant les principes du calcul 
diff\'erentiel d\'egag\'es de toute consid\'eration d'infiniment petits 
et d'\'evanouissans, de limites ou de fluxions 
et r\'eduits à l'analyse alg\'ebrique des quantit\'es finies 
}, Paris (an V). \\
Accessible on {\tt http://gallica.bnf.fr/}. 


\bibitem{EP} {\sc A. Lascoux}.
\emph{About the "$y$" in the $\chi_y$-characteristic of Hirzebruch}, 
  Conference in the honor of F. Hirzebruch, Institut Banach 1998,   
Contemp. Math. {\bf 241} (1999) 285--296.

\bibitem{Cbms} {\sc A. Lascoux}. 
\emph{Symmetric functions and combinatorial operators on polynomials},
CBMS/AMS Lectures Notes {\bf 99}, (2003).

\bibitem{Lothaire} {\sc A. Lascoux, B. Leclerc et J.Y. Thibon}.
\emph{The Plactic Monoid}, in 
\emph{Combinatorics on Words}, M. Lothaire ed, Cambridge Univ. Press (2002).

\bibitem{Littlewood} {\sc D.E. Littlewood}.
\emph{On certain symmetric functions}, Proc. London Math. Soc. 11 (1961)
485--498.

\bibitem{Mac} {\sc I.~G. Macdonald}.
\emph{Symmetric functions and {Hall} polynomials}.
Oxford University Press, Oxford, 2nd edition (1995).

\bibitem{Newton} {\sc I. Newton}. 
\emph{La m\'ethode des fluxions et des suites infinies,
par M. le Chevalier Newton}, 
pr\'eface et traduction de Buffon, 
Paris (1740).    

\bibitem{Warnaar}  {\sc O. Warnaar}. 
\emph{Hall-Littlewood functions and the A2 
 Rogers-Ramanujan  identities}, 
Adv. Math. {\bf 200} (2006) 403--434.


\end{thebibliography}
\end{document}